\documentstyle[12pt,amscd,xypic]{amsart}

\xyoption{all}
\CompileMatrices

\newcommand{\nc}{\newcommand}

\makeatletter
\@addtoreset{equation}{section}
\makeatother

\newtheorem{theorem}[equation]{Theorem}
\newtheorem{proposition}[equation]{Proposition}
\newtheorem{lemma}[equation]{Lemma}
\newtheorem{corollary}[equation]{Corollary}

\theoremstyle{definition}

\theoremstyle{remark}

\newtheorem{remark}[equation]{Remark}

\nc{\C}{{\mathbb{C}}}
\nc{\N}{{\mathbb{N}}}
\nc{\R}{{\mathbb{R}}}
\nc{\Z}{{\mathbb{Z}}}

\nc{\Mod}{{{\mathcal M}od}}

\nc{\FD}{{\mathfrak{D}}}
\nc{\FM}{{\mathfrak{M}}}
\nc{\CM}{{\mathcal{C}\!\mathcal{M}}}

\nc{\D}{{\mathcal{D}}}
\nc{\CE}{{\mathcal{E}}}
\nc{\CO}{{\mathcal{O}}}
\nc{\CU}{{\mathcal{U}}}
\nc{\CV}{{\underline{V}}}

\nc{\fg}{{\mathfrak{g}}}
\nc{\fgl}{{\mathfrak{gl}}}
\nc{\fu}{{\mathfrak{u}}}
\nc{\fpu}{{\mathfrak{pu}}}

\nc{\be}{{\mathbf{e}}}
\nc{\bn}{{\mathbf{n}}}
\nc{\bu}{{\mathbf{u}}}
\nc{\bv}{{\mathbf{v}}}
\nc{\bw}{{\mathbf{w}}}
\nc{\BK}{{\bar{K}}}
\nc{\BM}{{\mathbf{M}}}

\nc{\TG}{{\tilde{G}}}
\nc{\TZ}{{\tilde{Z}}}
\nc{\tx}{{\tilde{x}}}
\nc{\tbv}{{\tilde{\bv}}}
\nc{\tz}{{\tilde{\zeta}}}
\nc{\tmu}{{\tilde{\mu}}}

\nc{\urho}{\underline{\rho}}
\nc{\uB}{\underline{B}}
\nc{\ui}{\underline{i}}
\nc{\uj}{\underline{j}}

\nc{\eps}{\varepsilon}
\nc{\hrho}{{\hat{\rho}}}

\nc{\one}{{\mathbf{1}}}
\nc{\two}{{\mathbf{t}}}

\nc{\Rep}{{\mathop{\operatorname{\rm Rep}}}}
\nc{\Tot}{{\mathop{\operatorname{\rm Tot}}}}
\nc{\Ker}{{\mathop{\operatorname{\rm Ker}}}}
\nc{\Hilb}{{\mathop{\operatorname{\rm Hilb}}}}
\nc{\End}{{\mathop{\operatorname{\rm End}}}}
\nc{\Hom}{{\mathop{\operatorname{\rm Hom}}}}
\nc{\CHom}{{\mathop{\operatorname{{\mathcal{H}}\it om}}}}
\nc{\Coh}{{\mathop{\operatorname{\rm Coh}}}}
\nc{\GL}{{\mathop{\operatorname{\rm GL}}}}
\nc{\id}{{\mathop{\operatorname{\rm id}}}}
\nc{\rk}{{\mathop{\operatorname{\rm r}}}}
\nc{\supp}{{\mathop{\operatorname{\rm supp}}}}

\nc{\reg}{{\text{\rm reg}}}

\nc{\cplus}{{\mathbf{C}_+}}
\nc{\cminus}{{\mathbf{C}_-}}
\nc{\cthree}{{\mathbf{C}_*}}
\nc{\Qbar}{{\bar{Q}}}

\nc{\bh}{{\bar{h}}}
\nc{\bOmega}{{\overline{\Omega}}}

\nc{\seq}[1]{\stackrel{#1}{\sim}}

\title{Quiver varieties and Hilbert schemes}

\author{Alexander Kuznetsov}
\thanks{I was partially supported by RFFI grants 99-01-01144 and 99-01-01204.}
\address{
Institute for Problems of Information Transmission, 
 Russian Academy of Sciences,
19 Bolshoi Karetnyi, Moscow 101447, Russia
}
\email{sasha@@kuznetsov.mccme.ru}

\begin{document}

\begin{abstract}
In this note we give an explicit geometric description of some of the 
Nakajima's quiver varieties. More precisely, we show that the 
$\Gamma$-equivariant Hilbert scheme $X^{\Gamma[n]}$ and the Hilbert scheme 
$X_\Gamma^{[n]}$ (where $X=\C^2$, $\Gamma\subset SL(\C^2)$ is a finite 
subgroup, and $X_\Gamma$ is a minimal resolution of $X/\Gamma$) are quiver 
varieties for the affine Dynkin graph, corresponding to $\Gamma$ via 
the McKay correspondence, the same dimension vectors, but different 
parameters~$\zeta$ (for earlier results in this direction see \cite{H,VV,W}). 
In particular, it follows that the varieties $X^{\Gamma[n]}$ and 
$X_\Gamma^{[n]}$ are diffeomorphic. Computing their cohomology (in the 
case $\Gamma=\Z/d\Z$) via the fixed points of $(\C^*\times\C^*)$-action 
we deduce the following combinatorial identity: the number $UCY(n,d)$ 
of uniformly coloured in $d$ colours Young diagrams consisting of $nd$ 
boxes coincides with the number $CY(n,d)$ of collections of $d$ Young 
diagrams with the total number of boxes equal to $n$.
\end{abstract}

\maketitle

\section{Introduction}

The quiver varieties defined by Nakajima in \cite{Na} are of fundamental
importance in algebraic and differential geometry, theory of representations
and other branches of mathematics. They provide a rich source of examples
of hyperk\"ahler manifolds with very interesting geometry. For example,
the quiver variety $\FM_\zeta(\bv^0,0)$, corresponding to an affine 
Dynkin graph, the dimension vector $\bv^0$ given by the minimal positive 
imaginary root of the corresponding affine root lattice, 
and any generic parameter $\zeta$ is diffeomorphic to the 
minimal resolution $X_\Gamma$ of the simple singularity $X/\Gamma$, where
$X=\C^2$ and $\Gamma$ is a finite subgroup in $SL(\C^2)$, corresponding
to the graph via the McKay correspondence (see \cite{Kr}). Furthermore,
Nakajima proved that for some dimension vectors $\bv$, $\bw$ and 
for generic $\zeta$ the regular locus $\FM^{\text{reg}}_\zeta(\bv,\bw)$ 
of the quiver variety $\FM^{\text{reg}}_\zeta(\bv,\bw)$ is diffeomorphic 
to the framed moduli space of instantons on the 1-point compactification 
of $\FM_{-\zeta}(\bv^0,0)$. So, it would be natural to hope that one can
give an algebraic description of $\FM_\zeta(\bv,\bw)$ as some moduli space
of coherent sheaves on $\FM_{-\zeta}(\bv^0,0)$. The main result of this
note is a description of this kind in the case $\bv=n\bv^0$, 
$\bw=\bw^0$ --- the simple root of the extending vertex $0$ of the graph, 
and $\zeta=(0,\zeta_\R)$ with either $\zeta_\R\in\cplus$, 
or $\zeta_\R\in\cminus(n)$, where 
$$
\cplus = \{\zeta_\R\in\R^d\ |\ \zeta_\R^k>0\quad 0 \le k \le d-1\}
$$
$$
\cminus(n)=\left\{\zeta_\R\in\R^d\ \left|\ 
\frac1n\zeta_\R^k > \sum_{i=0}^{d-1}\zeta_\R^i\bv^0_i > 0
\quad 1\le k\le d-1\right.\right\}
$$
We prove that for $\zeta_\R\in\cplus$ the quiver variety 
$\FM_{(0,\zeta_\R)}(n\bv^0,\bw^0)$ is isomorphic to the 
$\Gamma$-equivariant Hilbert scheme of points on the plane~$X$.
This fact is well known, see e.g.~\cite{W}, \cite{VV}; we 
include a proof only for the sake of completeness.
We also prove that for $\zeta_\R\in\cminus(n)$ the quiver variety 
$\FM_{(0,\zeta_\R)}(n\bv^0,\bw^0)$ is isomorphic to 
the Hilbert scheme of points on $X_\Gamma$. After this paper
was written H.~Nakaijma kindly informed me that this fact 
was known to him (see~\cite{W} and~\cite{H}), but his arguments 
are different (see Remark~\ref{pn}). Our proof is based on 
an interpretation of quiver varieties as moduli spaces 
of representations of the corresponding double quivers 
suggested by Crawley-Boevey in~\cite{CB}.

The paper is organized as follows, In section~2 we recollect the necessary
background: the definition of quiver varieties, representations of quivers
and the construction of Crawley-Boevey. In section~3 we reproduce in 
a short form a geometric version of the McKay correspondence, based
on investigation of $X_\Gamma$. We also prove here a generalization 
of the result of Kapranov and Vasserot \cite{KV}. More precisely,
the authors of {\em loc.\ cit.} have constructed equivalences
of (bounded) derived categories
$$
\xymatrix{
\D^b(\Coh_\Gamma(X)) \ar[rr]<.5ex>^{\Phi} &&
\D^b(\Coh(X_\Gamma)) \ar[ll]<.5ex>^{\Psi}
},
$$
where $\Coh$ stands for the category of coherent sheaves,
and $\Coh_\Gamma$ stands for the category of $\Gamma$-equivariant
coherent sheaves. We show that there is a whole family of 
equivalences
$$
\xymatrix{
\D^b(\Coh_\Gamma(X)) \ar[rr]<.5ex>^{\Phi_{\zeta_\R}} &&
\D^b(\Coh(X_\Gamma)) \ar[ll]<.5ex>^{\Psi_{\zeta_\R}}
},
$$
which differ when parameters $\zeta_\R$ lie in distinct 
chambers of $\R^n$ with respect to roots hyperplanes.
In section~4 we prove the main results of this paper, isomorphisms
$$
\FM_{(0,\zeta_\R)}(n\bv^0,\bw^0)\cong 
\begin{cases} 
X^{\Gamma[n]},   & \text{for $\zeta_\R\in\cplus$}\\
X_\Gamma^{[n]}, & \text{for $\zeta_\R\in\cminus(n)$}
\end{cases}
$$
where $X^{\Gamma[n]}$ is the $\Gamma$-equivariant Hilbert scheme on $X$,
and $X_\Gamma^{[n]}$ is the Hilbert scheme on $X_\Gamma$. It follows
immediately from the general properties of quiver varieties that
$X^{\Gamma[n]}$ and $X_\Gamma^{[n]}$ are diffeomorphic.
We also prove a generalization of the first isomorphism in
the Calogero-Moser context: we check that for arbitrary $\tau\ne0$
and $\zeta_\R\in\cplus$ we have
$$
\FM_{((\tau,\dots,\tau),\zeta_\R)}(n\bv^0,\bw^0) \cong
(\CM_{nN})^\Gamma_{reg},
$$
where $N=|\Gamma|$, $\CM_{nN}$ is the Calogero-Moser space and
$(\CM_{nN})^\Gamma_{reg}$ is the connected component of the set
of $\Gamma$-fixed points, where the tautological bundle is 
a multiplicity of the regular representation.

Finally, in section~5 we consider $(\C^*\times\C^*)$-actions on
$X^{\Gamma[n]}$ and $X_\Gamma^{[n]}$ for cyclic $\Gamma\cong\Z/d\Z$.
We show that these actions have only finitely many number 
of fixed points and that the number of fixed points equal
to the dimension of the cohomology. Furthermore, we check that
the fixed points on $X^{\Gamma[n]}$ are in a bijection with the set of
uniformly coloured in $d$ colours Young diagrams with $dn$ boxes,
and that the fixed points on $X_\Gamma^{[n]}$ are in a bijection with
the set of collections of $d$ Young diagrams with total number of boxes
equal $n$. Recalling that $X^{\Gamma[n]}$ and $X_\Gamma^{[n]}$ are
diffeomorphic  we obtain the following combinatorial identity
(see \cite{JK} for a combinatorial proof)
$$
UCY(n,d) = CY(n,d),
$$
where $UCY$ and $CY$ denote the number of uniformly coloured diagrams
and the number of collections of diagrams respectfully.

\subsection*{Acknowledgements}
I am very grateful to V.~Ginzburg for bringing quiver varieties into
the area of my interests and to H.~Nakajima for helpful comments. 
Also I would like to thank M.~Finkelberg and D.~Kaledin for 
fruitful discussions.

\section{Quiver varieties}

In this section we recollect the definition and the basic properties 
of quiver varieties. We will follow the notation of \cite{Na}.

\subsection*{Definition}\label{defqv}

Choose an arbitrary finite graph and let $H$ be the
set of pairs, consisting of an edge together with an orientation on it.
Let $in(h)$ (resp.~$out(h)$) denote the incoming (resp.\ outgoing)
vertex of $h\in H$. Let $\bh$ denote the edge $h$ with the reverse 
orientation. Further, we choose an orientation of the graph, that 
is a subset $\Omega\subset H$ such that $\Omega\cup\bOmega=H$, 
$\Omega\cap\bOmega=\emptyset$. We define $\eps(h)=1$ for $h\in\Omega$
and $\eps(h)=-1$ for $h\in\bOmega$. Finally, we identify the set of 
vertices of our graph with the set $\{0,1,\dots,d-1\}$.

Choose a pair of hermitian vector spaces $V_k$, $W_k$ for each 
vertex of the graph and let 
$$
\bv=(\dim_\C V_0,\dots,\dim_\C V_{d-1}),\
\bw=(\dim_\C W_0,\dots,\dim_\C W_{d-1})\in\Z^d
$$
be their dimension vectors. Then the complex vector space
\begin{multline*}
\BM = \BM(\bv,\bw) = \\
\left(\bigoplus_{h\in H}\Hom(V_{out(h)},V_{in(h)})\right)\bigoplus
\left(\bigoplus_{k=0}^{d-1}\big[\Hom(W_k,V_k)\oplus\Hom(V_k,W_k)\big]\right) 
\end{multline*}
can be identified with the cotangent bundle of the hermitian vector space
$$
\BM_\Omega(\bv,\bw)=
\left(\bigoplus_{h\in\Omega}\Hom(V_{out(h)},V_{in(h)})\right)\oplus
\left(\bigoplus_{k=0}^{d-1}\Hom(W_k,V_k)\right).
$$
In particular, $\BM$ can be considered as a flat hyper-K\"ahler manifold.

Note that the group $G_\bv=\prod_{k=0}^{d-1} U(V_k)$ acts on $\BM$:
$$
g = (g_k)_{k=0}^{n-1}:
(B_h,i_k,j_k) \mapsto (g_{in(h)}B_hg_{out(h)}^{-1},g_ki_k,j_kg_k^{-1}),
$$
where $B_h\in\Hom(V_{out(h)},V_{in(h)})$, $i_k\in\Hom(W_k,V_k)$,
$j_k\in\Hom(V_k,W_k)$. This action evidently preserves the 
hyper-K\"ahler structure. The corresponding moment map 
$\mu=(\mu_\R,\mu_\C)$ is given by the following explicit formulas
$$
\mu_\C(B,i,j) = \left(\sum_{in(h)=k}\eps(h)B_hB_\bh+i_kj_k\right)
\in \bigoplus_{k=0}^{d-1}\fgl(V_k) = \fg_\bv\otimes\C,
$$
$$
\mu_\R(B,i,j) = 
\frac i2\left(\sum_{in(h)=k}
B_hB_h^\dagger-B_\bh^\dagger B_\bh+i_ki_k^\dagger-j_k^\dagger j_k
\right)\!\!\!\! \in \bigoplus_{k=0}^{d-1}\fu(V_k) = \fg_\bv,
$$
where $\fg_\bv$ is the Lie algebra of $G_\bv$ which is identified 
with its dual space $\fg_\bv^*$ via the canonical hermitian inner product.

Let $Z_\bv\subset\fg_\bv$ denote the center of the Lie algebra $\fg_\bv$. 
For any element $\zeta=(\zeta_\C,\zeta_\C)\in (Z_\bv\otimes\C)\oplus Z_\bv$ 
the corresponding quiver variety $\FM_\zeta$ is defined as 
a hyper-K\"ahler quotient
$$
\FM_\zeta = \FM_\zeta(\bv,\bw) = \{(B,i,j)\ |\ \mu(B,i,j)=-\zeta\}/G_\bv.
$$
In general, $\FM_\zeta$ has singularities, however its open subset
$$
\FM_\zeta^\reg = \left.\left\{(B,i,j)\in\mu^{-1}(-\zeta)\ \left|\ 
\parbox{0.34\textwidth}{the stabilizer of $(B,i,j)$ in $G_\bv$ is trivial}
\right.\right\}\right/G_\bv.
$$
is a smooth hyper-K\"ahler manifold (but maybe empty).

\subsection*{Roots and genericity}

The center $Z_\bv$ of $\fg_\bv=\bigoplus_{k=0}^{d-1}\fu(V_k)$ is the product 
of the set of scalar matrices on $V_k$, thus it can be considered as 
a subspace of $\R^d$.

Let $A$ denote the adjacency matrix of the graph and let $C=2I-A$
be the generalized Cartan matrix. Then we consider the set of 
positive roots
$$
R_+ = \{\theta = (\theta_k)\in(Z_{\ge0})^d\ |\ 
{}^t\theta C\theta\le 2\}\setminus\{0\}.
$$
Further we denote 
$$
R_+(\bv) = \{\theta\in R_+\ |\ 
\theta_k\le\bv_k\ \text{for all $k=0,\dots,d-1$}\},
$$
and for any positive root $\theta$ we consider a hyperplane
$$
D_\theta=\{x=(x_k)\in\R^d\ |\ \sum x_k\theta_k=0\}\subset \R^d.
$$
The element $\zeta$ is called generic (with respect to $\bv$) 
if for any $\theta\in R_+(\bv)$ we have 
$$
\zeta\not\in\R^3\otimes D_\theta\subset\R^3\otimes\R^d \supset
\R^3\otimes Z_\bv = (Z_\bv\otimes\C)\oplus Z_\bv.
$$
The importance of generic parameters $\zeta$ is explained by the
following Theorem of Nakajima (see \cite{Na} and \cite{CB} for 
connectedness).

\begin{theorem}\label{na}
For any $\zeta$, $\bv$, and $\bw$ the quiver variety $\FM_\zeta(\bv,\bw)$
is either empty, or nonemty and connected. Further, if $\zeta$ is generic 
then $\FM_\zeta$ is smooth and there is a canonical map
$\pi_0:\FM_{(\zeta_\C,\zeta_\R)}\to\FM_{(\zeta_\C,0)}$
which is a resolution of singularities, provided $\FM_{(\zeta_\C,0)}^\reg$
is nonempty. Finally, if both $\zeta$ and $\zeta'$ are generic then 
the varieties $\FM_\zeta(\bv,\bw)$ and $\FM_{\zeta'}(\bv,\bw)$ 
are diffeomorphic.
\end{theorem}

\subsection*{Representations of quivers}

Recall that a quiver $Q$ is a finite oriented graph.
Let $I(Q)$ denote the set of vertices and $A(Q)$ denote
the set of arrows of $Q$. For any arrow of the quiver 
$\alpha\in A(Q)$ we denote by $in(\alpha)$ and $out(\alpha)$ 
the incoming and outgoing vertex of $\alpha$. A representation 
$\rho$ of a quiver $Q$ is the following datum:
$$
\rho = (V,B),
$$
where $V = (V_i)_{i\in I(Q)}$ is a collection of vector spaces
for each vertex of the quiver, and $B = (B_\alpha)_{\alpha\in A(Q)}$
is a collection of linear maps $B_\alpha:V_{out(\alpha)}\to V_{in(\alpha)}$
for each arrow of the quiver. A morphism of representations
$$
\phi:\rho=(V,B)\to\rho'=(V',B')
$$
is a collection of linear maps $\phi_i:V_i\to V'_i$ for each vertex 
of the quiver, such that for any arrow $\alpha$ we have
$$
\phi_{in(\alpha)}B_\alpha = B'_\alpha\phi_{out(\alpha)}.
$$
The dimension of a representation $\rho=(V,B)$ is the collection
of dimensions of the vector space $V_i$:
$$
\dim\rho = (\dim V_i)_{i\in I(Q)}.
$$
For any dimension vector $\bv=(\bv_i)_{i\in I(Q)}$ let
$$
\Rep_Q(\bv) = \mathop{\bigoplus}_{\alpha\in A(Q)}
\Hom(\C^{\bv_{out(\alpha)}},\C^{\bv_{in(\alpha)}})
$$
the space of $\bv$-dimensional representations of $Q$.
The group 
$$
GL(\bv,\C) = \prod_{i\in I(Q)} GL(\bv_i,\C)
$$
acts on $\Rep_Q(\bv)$ by conjugation:
$$
g = (g_i)_{i\in I(Q)}:(B_\alpha)_{\alpha\in A(Q)} \mapsto
(g_{in(\alpha)}B_\alpha g_{out(\alpha)}^{-1}).
$$ 
This action evidently factors through the quotient group
$$
PGL(\bv,\C) = GL(\bv,\C)/\C^*,\qquad
$$
where the embedding $\C^*\to GL(\bv,\C)$ takes $\lambda\in\C^*$
to the element $(\mathop{\rm diag}(\lambda))_{i\in I(Q)}$.
Moreover, it is clear that the set of $PGL(\bv,\C)$ orbits
in $\Rep_Q(\bv)$ is the set of isomorphism classes of 
$\bv$-dimensional representations of $Q$.

Let $\chi$ be a map $I(Q)\to\R$, $\chi(i)=\chi^i$ (so-called polarization). 
For any dimension vector $\bv$ we define 
$$
\chi(\bv) = \sum_{i\in I(Q)}\chi^i\bv_i.
$$
A representation $\rho$ of a quiver $Q$ is called $\chi$-stable 
(resp.\ $\chi$-semistable) if $\chi(\dim\rho)=0$ and for any 
subrepresentation $\rho'\subset\rho$ such that $0\ne\rho'\ne\rho$ 
we have $\chi(\dim\rho') > 0$ (resp.\ $\chi(\dim\rho')\ge 0$). 
Representations $\rho$ and $\rho'$ are called $S$-equivalent with 
respect to a polarization $\chi$ if both $\rho$ and $\rho'$ are 
$\chi$-semistable and admit filtrations
$$
0 = \rho_0 \subset \rho_1 \subset \dots \subset \rho_n = \rho
\quad\text{and}\quad
0 = \rho'_0 \subset \rho'_1 \subset \dots \subset \rho'_n = \rho'
$$
such that $\chi(\dim\rho_i)=\chi(\dim\rho'_i)=0$ for all $i$, and
$$
\bigoplus_{i=1}^n\rho_i/\rho_{i-1} \cong
\bigoplus_{i=1}^n\rho'_i/\rho'_{i-1}.
$$

The following Theorem has been proved in \cite{Ki} and \cite{Na2}.

\begin{theorem}\label{mq}
For any quiver $Q$, dimension vector $\bv$ and polarization $\chi$,
such that $\chi(\bv)=0$, there exists a coarse moduli space 
$\Mod_Q(\bv,\chi)$ of $\bv$-dimensional $\chi$-semistable representations
of $Q$. Furthermore, if every $\chi$-semistable $\bv$-dimensional 
representation is $\chi$-stable and the dimension vector $\bv$ is 
indivisible, then $\Mod_Q(\bv,\chi)$ is a fine moduli space.
\end{theorem}

\subsection*{Quivers with relations}

Recall that the path algebra $\C[Q]$ of a quiver $Q$ is an algebra 
with a basis given by (oriented) paths in quiver and with multiplication 
given by concatenation of paths. It is clear that a representation 
of a quiver $\rho=(V,B)$ is the same as the structure of a right
$\C[Q]$-module on the vector space $\bigoplus_{i\in I(Q)}V_i$.

A quiver with relations is a pair $(Q,J)$, where
$Q$ is a quiver and $J$ is a two-sided ideal $J\subset\C[Q]$ 
in its path algebra. A representation of a quiver with relations 
$(Q,J)$ is a representation $\rho$ of its underlying quiver $Q$,
such that the ideal $J$ acts by zero in the corresponding 
right $\C[Q]$-module. In the other words, it is just a right
$\C[Q]/J$-module. We denote by $\Rep_{Q,J}(\bv)$ the space of 
$\bv$-dimensional representations of the quiver 
with relations $(Q,J)$. It is clear that $\Rep_{Q,J}(\bv)$ 
is a closed algebraic subset in the vector space $\Rep_Q(\bv)$.

It is easy to see that an analogue of Theorem~\ref{mq}
is true for quivers with relations.

\begin{theorem}\label{mqr}
For any quiver with relations $(Q,J)$, dimension vector $\bv$ and 
polarization $\chi$, such that $\chi(\bv)=0$, there exists a coarse 
moduli space $\Mod_{Q,J}(\bv,\chi)$ of $\bv$-dimensional $\chi$-semistable 
representations of $(Q,J)$. Furthermore, if every $\chi$-semistable 
$\bv$-dimensional representation is $\chi$-stable and the dimension 
vector $\bv$ is indivisible, then $\Mod_{Q,J}(\bv,\chi)$ is a fine 
moduli space.
\end{theorem}

\begin{remark}\label{qcham}
Let us fix the quiver $(Q,J)$ and the dimension vector $\bv$ and 
let $\chi$ vary in the space $D_\bv = \{\chi\ |\ \chi(\bv) = 0\}$.
Then it is easy to see that the notion of $\chi$-(semi)stability
can change only when the polarization $\chi$ crosses a hyperplane
$D_{\bv'}\cap D_{\bv} \subset D_{\bv}$ for some \hbox{$0 < \bv' < \bv$}.
Thus the space $D_\bv$ decomposes by walls $D_{\bv'}\cap D_{\bv}$
into a union of chambers. 
\end{remark}

\begin{remark}\label{explconstr}
One of constructions of the moduli space $\Mod_{Q,J}(\bv,\chi)$ is 
given in the terms of symplectic reduction. Namely, we can consider
the representations space $\Rep_{Q,J}(\bv)$ as a symplectic variety
with respect to the K\"ahler form (corresponding to the flat K\"ahler
metric on the vector space $\Rep_Q(\bv)$). Then compact algebraic group
$$
PU(\bv) = U(\bv)/U(1) = \left(\prod_{i\in I(Q)}U(\bv_i)\right)/U(1)
$$
(a real form of $PGL(\bv,\C)$) acts on $\Rep_{Q,J}(\bv)$ preserving
the symplectic structure. Let $\mu:\Rep_{Q,J}(\bv)\to\fpu(\bv)^*$
denote the corresponding moment map. The polarization $\chi$ can 
be thought of as an element of the center of dual Lie algebra $\fpu(\bv)^*$.
Then one can check (see e.g.~\cite{Ki}) that
$$
\Mod_{Q,J}(\bv,\chi) = \mu^{-1}(-\chi)/PU(\bv).
$$
\end{remark}

\subsection*{Modular description of quiver varieties}


The following construction, relating quiver varieties to
the moduli spaces of representation of quivers was suggested
in \cite{CB}. Assume that we are interested in a quiver variety
$\FM_\zeta(\bv,\bw)$ corresponding to a graph with the set of
vertices $\{0,\dots,d-1\}$ and with the set of oriented edges~$H$.
We define a quiver $Q = Q(\bw)$ (depending on the dimension 
vector $\bw$) with the set of vertices 
$$
I(Q) = \{*\}\cup\{0,\dots,d-1\}
$$
and with the following set of arrows. Firstly, we consider
every oriented edge $h\in H$ as an arrow in the quiver with
the same outgoing and incoming vertices. Secondly, for every 
$0\le i\le d-1$ we draw $\bw_i$ arrows from the vertex $*$ 
to the vertex $i$ and backwards.

Then it is easy to see that a choice of a basis in spaces 
$W_0$, \dots, $W_{d-1}$ gives an identification
$$
\BM(\bv,\bw) = \Rep_{Q(\bw)}((1,\bv)).
$$
Further, the expression $\mu_\C + \zeta_\C$ can 
be considered as an element of the path algebra 
$\C[Q(\bw)]$ of the quiver~$Q(\bw)$. We denote 
by $J(\zeta_\C)$ the two-sided ideal in $\C[Q(\bw)]$, 
generated by $\mu_\C+\zeta_\C$. Then the algebraic
subvariety $\mu_\C^{-1}(-\zeta_\C)\subset\BM(\bv,\bw)$
gets identified with the algebraic subvariety
$\Rep_{Q(\bw),J(\zeta_\C)}((1,\bv)) \subset \Rep_{Q(\bw)}((1,\bv))$:
\begin{equation}\label{mci}
\mu_\C^{-1}(-\zeta_\C) = \Rep_{Q(\bw),J(\zeta_\C)}((1,\bv)).
\end{equation}
Further, note that the quiver variety
$$
\FM_\zeta(\bv,\bw) = \mu_\C^{-1}(-\zeta_\C)\cap\mu_\R^{-1}(-\zeta_\R)/G_\bv
$$
is nothing but the symplectic reduction of the variety
$\mu_\C^{-1}(-\zeta_\C) = \Rep_{Q(\bw),J(\zeta_\C)}((1,\bv))$
with respect to its canonical K\"ahler form. 
Finally, the group $G_\bv$ acting on $\mu_\C^{-1}(-\zeta_\C)$ 
is canonically isomorphic to the group $PU((1,\bv))$ acting
on $\Rep_{Q(\bw),J(\zeta_\C)}((1,\bv))$:
$$
G_{\bv} = \prod_{i=0}^{d-1}U(\bv_i) \cong 
\left(U(1)\times\prod_{i=0}^{d-1}U(\bv_i)\right)/U(1),
$$
and under this identificationthe element
$\zeta_\R\in\fg_\bv$ goes to the polarization
\begin{equation}\label{chiz}
\chi_\bv(\zeta_\R) = (-\zeta_\R(\bv),\zeta_\R)\in \fpu^*((1,\bv)).
\end{equation}
Summing up and taking into account Remark~\ref{explconstr} we obtain 
the following Proposition.

\begin{proposition}[cf.\ \cite{CB}]\label{md}
For any $(\zeta_\C,\zeta_\R)$ we have an isomorphism 
of algebraic varieties
\begin{equation}\label{id}
\FM_{(\zeta_\C,\zeta_\R)}(\bv,\bw) = 
\Mod_{Q(\bw),J(\zeta_\C)}((1,\bv),\chi_\bv(\zeta_\R)).
\end{equation}
\end{proposition}

\begin{remark}\label{qvcham}
The arguments of Remark~\ref{qcham} show that $\FM_\zeta(\bv,\bw)$
depend on $\zeta_\R$ as follows. The space $\R^n$ of all $\zeta_\R$
is decomposed by root hyperplanes (walls) into a union of chambers.
Whenever $\zeta_\R$ varies within a chamber the complex structure
of $\FM_{(\zeta_\C,\zeta_\R)}(\bv,\bw)$ doesn't change at all
(but the hyper-K\"ahler metric does), while when $\zeta_\R$ crosses 
a wall $\FM_{(\zeta_\C,\zeta_\R)}(\bv,\bw)$ endures a (usually 
birational) transformation.
\end{remark}

We will finish this section with the following useful Lemma.

\begin{lemma}\label{fine}
If $\zeta$ is generic then 
$\Mod_{Q(\bw),J(\zeta_\C)}((1,\bv),\chi_\bv(\zeta_\R))$
is a fine moduli space.
\end{lemma}
{\sl Proof:} 
Let $\chi = \chi_\bv(\zeta_\R)$.
Since the dimension vector $(1,\bv)$ is evidently indivisible
it suffices to check that every $\chi$-semistable $(1,\bv)$-dimensional 
representation $\rho$ is $\chi$-stable. So assume that $\rho'\subset\rho$ 
is a $\chi$-stable subrepresentation, such that $0\ne\rho'\ne\rho$ and 
$\chi(\dim\rho')=0$. We have two possibilities: either $\dim\rho'=(0,\bv')$ 
or $\dim\rho'=(1,\bv')$ for some $\bv'$. 

Assume for example that $\dim\rho'=(0,\bv')$. Then we have
$$
0 = \chi(\dim\rho') =
(-\zeta_\R(\bv),\zeta_\R)(0,\bv') = \zeta_\R(\bv').
$$
Without loss of generality we can assume that $\rho'$ is $\chi$-stable.
It follows from Theorem~0.2 of \cite{CB} that $(0,\bv')$ is a positive 
root for the quiver $Q(\bw)$, hence $\bv'\in R_+(\bv)$, hence 
$\zeta_\R\in D_{\bv'}$. On the other hand, by Proposition~\ref{md}
the representation $\rho'$ corresponds to a point $(B',0,0)$ 
of the quiver variety $\FM_\zeta(\bv',\bw)$, hence
\begin{multline*}
\zeta_\C(\bv') = \sum_{k=0}^{n-1} Tr(\zeta_\C^k) = 
\sum_{k=0}^{n-1} Tr(-\mu_\C(B',0,0)) = \\ =
- Tr\left(\sum_{h\in H}\eps(h)B'_hB'_\bh\right) =
- Tr\left(\sum_{h\in \Omega}\big[B'_h,B'_\bh\big]\right) = 0.
\end{multline*}
Thus $\zeta\in\R^3\otimes D_{\bv'}$, a contradiction with 
the assumtion that $\zeta$ is generic.

Similarly, assume that $\dim\rho'=(1,\bv')$.
Consider the quotient representation $\rho'' = \rho/\rho'$.
Then $\dim\rho'' = (0,\bv'')$ with $\bv'' = \bv-\bv'$ and
without loss of generality we can assume that $\rho''$ is $\chi$-stable.
Then repeating the above arguments for $\rho''$ and $\bv''$ instead
of $\rho'$ and $\bv'$ we again deduce that $\zeta$ is not generic.

Thus we have proved that every $\chi$-semistable representation is 
$\chi$-stable, hence the moduli space is fine.
\qed

\section{The McKay correspondence}

The McKay correspondence associates to every (conjugacy class of~a) 
finite subgroup of the group $SL(\C^2)$ an affine Dynkin graph of 
type $ADE$. The graph corresponding to a subgroup $\Gamma\subset SL(\C^2)$
can be constructed as follows. Let $R_0$, \dots, $R_{d-1}$ be the set of all 
(isomorphism classes of) irreducible representations of the group $\Gamma$, 
and assume that $R_0$ is the trivial representation. Let $L=\C^2$ be 
the tautological representation. Further, let $a_{k,l}$ be the 
multiplicities in the decomposition of the tensor product 
$R_k\otimes L$ into the sum of irreducible representations:
\begin{equation}\label{rl}
R_k\otimes L \cong \mathop{\bigoplus}_{l=0}^{d-1} R_l^{\oplus a_{k,l}}.
\end{equation}
Then the graph with the set of vertices $\{0,\dots,d-1\}$
and with $a_{k,l}$ edges between the vertices $k$ and $l$
is the corresponding affine Dynkin graph. 

\subsection*{Quiver varieties for affine Dynkin graphs}

From now on we will consider quiver varieties corresponding
to an affine Dynkin graph. Let $\Gamma$ be the corresponding 
subgroup of $SL(\C^2)$. Let $d$ denote the number of irreducible
representations of $\Gamma$ (that is equal to the number
of vertices in the graph), and let $N$ be the order of 
the group $\Gamma$.

Let $V$ and $W$ be representations of $\Gamma$. Then
\begin{equation}\label{uudecomp}
V = \bigoplus_{k=0}^{d-1}V_k\otimes R_k,\qquad
W = \bigoplus_{k=0}^{d-1}W_k\otimes R_k.
\end{equation}
Then the dimension vectors
$$
\bv = (\dim V_0,\dots,\dim V_{d-1}),\quad
\bw = (\dim W_0,\dots,\dim W_{d-1}).
$$
can be thought of as the classes of $V$ and $W$ in 
the Grothendieck ring $K_0(\Gamma)$. Consider a triple
$$
(B,i,j) \in \Hom_\Gamma(V\otimes L,V) \oplus 
            \Hom_\Gamma(W,V) \oplus \Hom_\Gamma(V,W)
$$
where $\Hom_\Gamma$ denotes the space of all $\Gamma$-equivariant 
linear maps. Then from~(\ref{uudecomp}) and from
$$
V\otimes L = \bigoplus_{k=0}^{d-1}V_k\otimes R_k\otimes L =
\bigoplus_{k,l=0}^{d-1}V_k\otimes R_l^{\oplus a_{k,l}} =
\bigoplus_{k,l=0}^{d-1}V_k^{\oplus a_{k,l}}\otimes R_l
$$
it follows that a choice of $(B,i,j)$ is equivalent to
a choice of a collection $(B_h,i_k,j_k)\in\BM(\bv,\bw)$.
Thus
$$
\BM(\bv,\bw) = \Hom_\Gamma(V\otimes L,V) \oplus 
               \Hom_\Gamma(W,V) \oplus \Hom_\Gamma(V,W).
$$
Further, it is easy to check that 
$$
\mu_\C(B,i,j) = [B,B] + ij,\qquad
\mu_\R(B,i,j) = B B^\dagger - B^\dagger B + i i^\dagger - j^\dagger j,
$$
where $[B,B]$, $B B^\dagger$, and $B^\dagger B$ are defined 
as the following compositions
$$
\begin{array}{rcl}
[B,B] &:& V @>{\lambda_L^{-1}}>> 
V\otimes L\otimes L @>{B\otimes1}>> V\otimes L @>{B}>> V,\\
B B^\dagger &:& V @>{h_V}>> V^* @>{B^*}>> V^*\otimes L^*
@>{h_V^{-1}\otimes h_L^{-1}}>> V\otimes L @>{B}>> V,\\
B^\dagger B &:& V @>{B}>> V\otimes L^* @>{h_V\otimes h_L^{-1}}>>
V^*\otimes L @>{B^*}>> V^* @>{h_V^{-1}}>> V,
\end{array}
$$
and $\lambda_L$, $h_L$ and $h_V$ stand for $\Gamma$-invariant symplectic 
form and hermitian inner products on $L$ and $V$ respectively. Further, 
we can consider the parameters $\zeta_\C$ and $\zeta_\R$ as elements 
of $Z(\C[\gamma])$  and $\Z(\R[\Gamma])$, the center of the group 
algebra of $\Gamma$ over $\C$ and $\R$ respectively (element $\zeta_\C$ 
acts in $R_k$ as $\zeta_\C^k$-multiplication). Finally, it is easy to see 
that $G_\bv = U_\Gamma(V)$ and we obtain the following Lemma.

\begin{lemma}\label{gqv}
We have
$$
\FM_{(\zeta_\C,\zeta_\R)}(\bv,\bw) = \left.\left\{(B,i,j)\ \left|\ 
\begin{array}{rcl}
{}[B,B]+ij &=& - \zeta_\C\\
{}[B,B^\dagger] + i i^\dagger - j^\dagger j &=& -\zeta_\R
\end{array}
\right.\right\}\right/U_\Gamma(V).
$$
\end{lemma}

Let $X=\C^2$ with tautological action of $\Gamma$ and
identify $L$ with the space of linear functions on $X$.

\begin{lemma}[cf.\ \cite{KV}, Pr.~3.4]\label{vvk}
For any point $(B,i,j)\in\FM_{(0,\zeta_\R)}(\bv,\bw)$ such that $j=0$,
the map $B$ induces the structure of a $\Gamma$-equivariant 
$\C[X]$-module on $V$ {\rm(}where the linear functions 
on $X$ act via $B${\rm)}.
\end{lemma}
{\sl Proof:}
We have
$$
[B,B] = \mu_\C(B,i,j) - ij = 0 - ij = 0,
$$
hence the actions of linear functions on $X$ commute. Therefore,
they induce an action of $\C[X]$. The obtainde $\C[X]$-module 
structure on $V$ is evidently $\Gamma$-equivariant.
\qed\medskip

Combining this with (\ref{mci}) we obtain the following.

\begin{corollary}\label{vvk1}
We have an isomorphism
\begin{multline*}
\{(V_\bullet,B,i,j)\in\Rep_{Q(\bw),J(0)}((1,\bv))\ |\ j=0 \} \cong \\
\left\{\parbox{0.85\textwidth}{\raggedright%
$\Gamma$-equivariant $\C[X]$-module structures 
on $V = \bigoplus\limits_{k=0}^{d-1}V_k\otimes R_k$\\ 
with a framing $i\in\Hom_\Gamma(W,V)$}\right\}
\end{multline*}
\end{corollary}

Abusing the notation we will denote by $V_\bullet$ the representation 
of the quiver $(Q(\bw),J(0))$ corresponding to a $\Gamma$-equivariant 
$\C[X]$-module $V$ and $\Gamma$-morphism $i:W\to V$. Vice versa,
the $\Gamma$-equivariant $\C[X]$-module, corresponding to a
representation $\rho$ of the quiver $(Q(\bw),J(0))$ will be denoted 
by~$\Tot(\rho)$.

\begin{remark}\label{morita}
In fact, one can reformulate Corollary~\ref{vvk1} in terms of 
Morita equivalence. Let $e_*\in\C[Q(\bw)]$ denote the idempotent
of the vertex $*$ (the path of length $0$ with $*$ being
the outgoing (and incoming) vertex). Then the algebra 
$\C[Q(\bw)]/\langle J(0),e_*\rangle$ is Morita-equivalent 
to the smash product algebra $\C[X]\#\Gamma$.
\end{remark}

\subsection*{Resolutions of simple singularities}

We fix a pair of dimension vectors:
$$
\bv^0 = (\dim R_0,\dim R_1,\dots,\dim R_{d-1}),\quad\text{and}\quad
\bw^0 = (1,0,\dots,0).
$$
Note that 
\begin{equation}\label{n}
N = \sum_{k=0}^{d-1}(\bv^0_i)^2 = |\Gamma|.
\end{equation}
We fix also the complex parameter $\zeta_\C=0$ but let 
the real parameter $\zeta_\R$ vary. The quiver varieties 
$\FM_{(0,\zeta_\R)}(\bv^0,\bw^0)$ are described by the following 
Theorem.

\begin{theorem}[see \cite{Kr}]\label{kr}
$(i)$ The quiver variety $\FM_{(0,0)}(\bv^0,\bw^0)$ is isomorphic
to the quotient variety $X/\Gamma$.      

\noindent$(ii)$ 
For any generic $\zeta_\R$ the quiver variety 
$X_\Gamma(\zeta_\R)=\FM_{(0,\zeta_\R)}(\bv^0,\bw^0)$ 
is a minimal resolution of the quotient variety 
$X/\Gamma$ via the canonical map 
$\pi_0:\FM_{(0,\zeta_\R)}(\bv^0,\bw^0)\to\FM_{(0,0)}(\bv^0,\bw^0)$.
\end{theorem}

\begin{remark}
Actually, results of \cite{Kr} concern varieties $\FM_\zeta(\bv^0,0)$
rather than $\FM_\zeta(\bv^0,\bw^0)$. However, it is easy to see that these
varieties are canonically isomorphic. Indeed, arguing like in the proof 
of Lemma~\ref{li0j0} below, one can check that for any point of 
$\FM_\zeta(\bv^0,\bw^0)$ we have 
$$
j_0i_0 = \sum\zeta_\C^k\bv^0_k,\qquad
i_0i_0^\dagger - j_0^\dagger j_0 = -2\sqrt{-1}\sum\zeta_\R^k\bv^0_k.
$$
It follows that the map
$\BM(\bv^0,\bw^0)\supset\mu^{-1}(\zeta)\to\mu^{-1}(\zeta)\subset\BM(\bv^0,0)$,
induced by the canonical projection $\BM(\bv^0,\bw^0)\to\BM(\bv^0,0)$
(forgeting of $i_0$ and $j_0$) is a principal $U(1)$-bundle. Moreover,
this map is $G_{\bv^0}$-equivariant and the action of $U(1)$ on fibers
is nothing but the action of a subgroup $U(1)\subset G_{\bv^0}$. Hence
the quotients modulo $G_{\bv^0}$ are canonically isomorphic.
\end{remark}

\begin{remark}
In dimension 2 any two minimal resolutions necessarily coincide.
Thus, we have a canonical identification 
$X_\Gamma(\zeta_\R)\cong X_\Gamma(\zeta'_\R)$ for any 
generic $\zeta_\R$ and $\zeta'_\R$, compatible with the
projection to $X/\Gamma$. So, we can (and will) write 
$X_\Gamma$ instead of $X_\Gamma(\zeta_\R)$ without
risk of misunderstanding.
\end{remark}

Theorem~\ref{md} implies that for any $\zeta_\R$ we have an isomorphism
$$
X_\Gamma = 
\Mod_{Q(\bw),J(0)}((1,\bv^0),\chi_{\bv^0}(\zeta_\R)),
$$
while Lemma~\ref{fine} implies that $X_\Gamma$ is a fine moduli space
whenever $\zeta_\R$ is generic. This means that for any generic $\zeta_\R$ 
we have a universal representation $\urho^{\zeta_\R}$ of the quiver 
$(Q(\bw^0),J(0))$ over $X_\Gamma$, that is a collection of vector 
bundles $\CV_k^{\zeta_\R}$ over $X_\Gamma$, 
$k\in\{*\}\cup\{0,\dots,d-1\}$, which ranks are given by
\begin{equation}\label{rvk}
\rk(\CV_*^{\zeta_\R})=1,\qquad
\rk(\CV_k^{\zeta_\R})=\bv^0_k\quad\text{for $k\in\{0,\dots,d-1\}$,}
\end{equation}
and morphisms 
$$
\uB_h^{\zeta_\R}:\CV_{out(h)}^{\zeta_\R}\to\CV_{in(h)}^{\zeta_\R},\quad
\ui_0^{\zeta_\R}:\CV_*^{\zeta_\R}\to\CV_0^{\zeta_\R},\quad
\uj_0^{\zeta_\R}:\CV_0^{\zeta_\R}\to\CV_*^{\zeta_\R}
$$
(recall the choice of $\bw^0$), satisfying the equations
\begin{equation}\label{mmc}
\begin{array}{rl}
\sum_{in(h)=k} \eps(h)\uB_h^{\zeta_\R}\uB_\bh^{\zeta_\R} = 0 & 
\text{for $k=1,\dots,d-1$} \smallskip\\
\ui_0^{\zeta_\R}\uj_0^{\zeta_\R} +
\sum_{in(h)=0} \eps(h)\uB_h^{\zeta_\R}\uB_\bh^{\zeta_\R} = 0, 
\end{array}
\end{equation}
and such that for any $x\in X_\Gamma$ the corresponding representation
$$
\urho(x)^{\zeta_\R} = 
({\CV}_{|x}^{\zeta_\R},\uB^{\zeta_\R}(x),
\ui_0^{\zeta_\R}(x),\uj_0^{\zeta_\R}(x))
$$ 
of the quiver $Q(\bw^0)$ is $\chi_{\bv^0}(\zeta_\R)$-stable. 

\begin{remark}\label{rcham}
In fact, the family $\rho^{\zeta_\R}$ doesn't
change whenever $\zeta_\R$ varies within a chamber.
\end{remark}

Note that universal representation is defined up to a twist 
by a line bundle, so without loss of generality we may (and will) 
assume that $\CV_*^{\zeta_\R}=\CO$. 

\begin{lemma}\label{li0j0}
For any $(1,\bv)$-dimensional representation $(V,B_h,i_0,j_0)$ of the quiver
$(Q(\bw^0),J(0))$ we have $j_0i_0 = 0$.
\end{lemma}
{\sl Proof:}
Since $\dim V_* = 1$ it follows that
$$
j_0i_0 = \mathop{Tr} j_0i_0 = \mathop{Tr} i_0j_0.
$$
On the other hand, summing up equations~(\ref{mmc}) and taking 
the trace we see that
$$
\mathop{Tr}i_0j_0 =
\sum_{h\in\Omega} \mathop{Tr}[B_h,B_\bh] + \mathop{Tr}i_0j_0 =
\mathop{Tr}\mu_\C(B_h,i_0,j_0) = 0.
$$
\qed\medskip

Since $\rk(\CV^{\zeta_\R}_*) = \rk(\CV^{\zeta_\R}_0) = 1$ by (\ref{rvk}), 
it follows from Lemma~\ref{li0j0} that for any $x\in X_\Gamma$ 
either we have $\uj_0^{\zeta_\R}(x)=0$, or $\ui_0^{\zeta_\R}(x)=0$. 
Assume that
\begin{equation}\label{pos}
\zeta_\R(\bv) > 0.
\end{equation}
Then $\ui_0^{\zeta_\R}(x)=0$ would violate 
the $\chi_{\bv^0}(\zeta_\R)$-stability, 
because then $\urho^{\zeta_\R}(x)$ would admit 
a $(1,0,\dots,0)$-dimensional subrepresentation, but
$$
\chi_{\bv^0}(\zeta_\R)(1,0,\dots,0) = 
(-\zeta_\R(\bv^0),\zeta_\R)(1,0,\dots,0) = -\zeta_\R(\bv^0)<0.
$$ 
So it follows that when (\ref{pos}) holds we have 
$\uj_0^{\zeta_\R}=0$, and that $\ui_0^{\zeta_\R}$ is an embedding 
of line bundles, and hence an isomorphism. Thus we have proved
the following Lemma.

\begin{lemma}\label{norm}
For any generic $\zeta_\R$ satisfying $(\ref{pos})$ there exists 
a universal $\chi_{\bv^0}(\zeta_\R)$-stable family on $X_\Gamma$ such that
$$
\CV_*^{\zeta_\R}\cong\CV_0^{\zeta_\R}\cong\CO_{X_\Gamma},\qquad
\uj_0^{\zeta_\R}=0,\quad\text{and}\quad \ui_0^{\zeta_\R}=\id.
$$
\end{lemma}

Applying Corollary~\ref{vvk1} we deduce.

\begin{corollary}
For any generic $\zeta_\R$ such that $(\ref{pos})$ holds,
$\Tot(\rho^{\zeta_\R})$ is a family of $\Gamma$-equivariant
$\C[X]$-modules over $X_\Gamma$.
\end{corollary}

From now on we will denote by $\urho^{\zeta_\R}$ the universal
representation normalized as in Lemma~\ref{norm}. Moreover, to 
unburden the notation we will omit the superscript $\zeta_\R$
when it won't lead us to a confusion.

Let $\rho_0$ denote the unique $(1,0,\dots,0)$-dimensional 
representation of the quiver $(Q(\bw^0),J(\zeta_\C))$. It follows
from Lemma~\ref{norm} that for any $x$ we have canonical surjective
homomorphism $\urho(x)\to\rho_0$. Let us denote its kernel by
$\hrho(x)$. 

\subsection*{Derived equivalences}

For each generic $\zeta_\R$ we define functors
$$
\xymatrix{
\D^b(\Coh_\Gamma(X)) \ar[rr]<.5ex>^{\Phi_{\zeta_\R}} &&
\D^b(\Coh(X_\Gamma)) \ar[ll]<.5ex>^{\Psi_{\zeta_\R}}
},
$$
where 
$\Coh_\Gamma(X)$ denotes the category of $\Gamma$-equivariant 
coherent sheaves on $X$, 
$\Coh(X_\Gamma)$ denotes the category of coherenet sheaves on $X_\Gamma$, 
and $\D^b$ stands for the bounded derived category.
To this end we consider the family $\Tot(\urho^{\zeta_\R})$ of 
$\Gamma$-equivariant $\C[X]$-modules as a sheaf on the product 
$X_\Gamma\times X$. Denoting the projections to $X_\Gamma$ and 
to $X$ by $p_1$ and $p_2$ respectfully, we define
$$
\begin{array}{l}
\Phi_{\zeta_\R}(F) = (Rp_{1*}R\CHom(\Tot(\urho^{\zeta_\R}),p_2^*F))^G,
\quad\text{and}\\
\Psi_{\zeta_\R}(F) = Rp_{2*}(p_1^*F\otimes^L\Tot(\urho^{\zeta_\R}))
\end{array}
$$

\begin{theorem}
For any generic $\zeta_\R$ the functors $\Phi_{\zeta_\R}$ and 
$\Psi_{\zeta_\R}$ are mutually inverse equivalences of categories.
\end{theorem}

This Theorem is nothing but a slight generalization of Theorem~1.4
of Kapranov and Vasserot (see \cite{KV}). In fact, the equivalences
$\Phi$ and $\Psi$ from {\em loc.\ cit.} are isomorphic to the functors
$\Phi_{\zeta_\R}$ and $\Psi_{\zeta_\R}$ with 
$$
\zeta_\R\in \{ \zeta_\R\in \R^n\ |\ \zeta_\R^i<0
\text{ for all $0\le i\le n-1$}\}
$$

The proof of this Theorem can be done by the same arguments as 
in {\em loc.\ cit.}

\begin{remark}
In fact, the equivalences $\Phi_{\zeta_\R}$ and $\Psi_{\zeta_\R}$
don't change whenever $\zeta_\R$ varies within a chamber 
(see Remark~\ref{rcham}).
\end{remark}

\begin{remark}
For any pair of generic $\zeta_\R$ and $\zeta'_\R$ the compositions 
of functors $\Phi_{\zeta_\R}\cdot\Psi_{\zeta'_\R}$ and
$\Psi_{\zeta'_\R}\cdot\Phi_{\zeta_\R}$ are autoequivalences 
of the derived categorie $\D^b(\Coh(X_\Gamma))$ and $\D^b(\Coh_\Gamma(X))$
respectively. One can check that they generate the action of the affine
braid group described in~\cite{ST}.
\end{remark}

\section{Interpretation of quiver varieties}

\subsection*{Symmetric power of $X/\Gamma$}

From now on we will be interesting in an explicit
geometric description of quiver varieties $\FM_\zeta(n\bv^0,\bw^0)$
for $\zeta_\C=0$, $\bv^0$ and $\bw^0$ as above and various
$n$ and $\zeta_\R$. We begin with the simplest case,
$\zeta_\R=0$. In fact, this is well known, 
but we put here a proof for the sake of completeness.

\begin{proposition}\label{psi0}
We have $\FM_{(0,0)}(n\bv^0,\bw^0) = S^n(X/\Gamma)$,
where $S^n$ stands for the symmetric power.
\end{proposition}
{\sl Proof:}
The case $n=1$ follows from the first assertion of Theorem~\ref{kr}. 
So assume that $n>1$.

It follows from Proposition~\ref{md} that it suffices to check that
$$
\Mod_{Q(\bw),J(0)}((1,n\bv^0),0) = S^n(X/\Gamma).
$$
Fix some generic $\zeta_\R$. For a collection of points 
$(x_1,\dots,x_n)$, $x_i\in X/\Gamma$ we put
$$
g_0(x_1,\dots,x_n) = \rho_0 \oplus\hrho(\tx_1)\oplus\dots\hrho(\tx_n),
$$
where $\tx_i\in\pi_0^{-1}(x_i)\in X_\Gamma(\zeta_\R)$ are arbitrary lifts
of the points $x_i$. Since every representation is semistable 
with respect to the trivial polarization $\chi=0$, it follows 
that $g_0$ induces a map
$$
g_0:S^n(X/\Gamma) \to \Mod_{Q(\bw^0),J(0)}((1,n\bv^0),0).
$$
To construct the inverse we need the following Lemma.

\begin{lemma}\label{s0}
Any $(1,\bv)$-dimensional representation of $(Q(\bw^0),J(0))$
is $S$-equivalent {\rm(}with respect to $\chi=0${\rm)} 
to a representation with $j=0$.
\end{lemma}
{\sl Proof:}
Let $\rho=(V_*,V_\bullet)$ be a $(1,\bv)$-dimensional representation 
of $(Q(\bw^0),J(0))$. Let $B:V\otimes L\to V$, $i:W\to V$,
and $j:V\to W$ be the corresponding $\Gamma$-equivariant morphisms. 
Let $U$ be the minimal subspace of $V$ such that $i(W)\subset U$ 
and $B(U\otimes L)\subset U$. Then $U$ is invariant under the action 
of $\Gamma$, hence $U=\oplus U_l\otimes R_k$. Moreover, 
$(V_*,U_\bullet)$ is a subrepresentation of $(V_*,V_\bullet)$, thus
$$
(V_*,V_\bullet) \seq0 (V_*,U_\bullet) \oplus (0,V_\bullet/U_\bullet),
$$
where $\seq0$ denotes $S$-equivalence with respect to $\chi=0$.
On the other hand, the arguments of \cite{Na2}, Lemma~2.8 show
that $j=0$ in the first summand in the RHS, and in the second
summand $j=0$ for trivial reasons. 
\qed\medskip

Now we need to construct a map inverse to $g_0$. 
Take arbitrary $\rho\in\Mod_{Q(\bw^0),J(0)}((1,n\bv^0),0)$.
By Lemma~\ref{s0} we can assume that $j=0$ in $\rho$.
Then $\Tot(\rho)$ is a $\Gamma$-equivariant $\C[X]$-module.
Regarding it as a $\Gamma$-equivariant coherent sheaf on $X$
we define
$$
f_0(\rho) = \supp(\Tot(\rho)) \in (S^{nN}X)^\Gamma = S^n(X/\Gamma).
$$
Note that the definition of the map $f_0$ is correct. Indeed, assume
that $\rho'\seq0\rho$ and $j=0$ in $\rho'$. Then $\rho$ and $\rho'$
admit filtrations with isomorphic associated factors (up to a permutation).
These filtrations induce filtrations of sheaves $\Tot(\rho)$ and
$\Tot(\rho')$ with isomorphic associated factors. But then
$\supp(\Tot(\rho)) = \supp(\Tot(\rho'))$.

Now we have to check that the maps $g_0$ and $f_0$ are 
mutually inverse. More precisely, one have to check that
\begin{equation}\label{todo}
g_0(f_0(\rho)) \seq0 \rho\quad\text{and}\quad
f_0(g_0(x_1,\dots,x_n)) = (x_1,\dots,x_n).
\end{equation}
Before we begin the proof let us note that from the explicit 
construction of isomorphisms in Theorem~\ref{kr} it follows
that for any point $\tx\in X_\Gamma$ we have
$$
\supp\Tot(\hrho(\tx)) = \pi_0(\tx).
$$

Now, take $\rho\in\Mod_{Q(\bw^0),J(0)}((1,n\bv^0),0)$.
By Lemma~\ref{s0} we can assume that $j=0$ in $\rho$,
hence $\rho\seq0\rho_0\oplus\hrho$. Now consider the 
sheaf $\Tot(\hrho)$. It is a $\Gamma$-equivariant sheaf
of finite length on $X$, hence it admits a filtration 
with associated factors being either the structure
sheaves of a $\Gamma$-orbit corresponding to a point 
$0\ne x_i\in X/\Gamma$, or a length 1 skyscraper sheaf 
supported at zero. This filtration induces a filtration 
on $\hrho$ with associated factors being either $\hrho(\tx_i)$, with 
$\tx_i = \pi_0^{-1}(x_i)$, or $\rho_k$, 
where $\rho_k$ is the unique representation
with $\dim(\rho_k)=\bw^k = (0,\dots,0,1,0,\dots,0)$ 
($1$ stands on the $k$-th position). Thus
\begin{equation}\label{seq0}
\rho \seq0 
\rho_0 \oplus 
\left(\bigoplus_{i=1}^m\hrho(\tx_i)\right) \oplus
\left(\bigoplus_{k=0}^{d-1}(\rho_k)\right)^{\oplus(n-m)\bv^0_k}
\end{equation}
Then
\begin{multline*}
g_0(f_0(\rho)) = g_0(\supp\Tot(\rho)) = 
g_0(x_1,\dots,x_m,0,\dots,0) = \\
\rho_0 \oplus 
\left(\bigoplus_{i=1}^m\hrho(\tx_i)\right) \oplus
\hrho(\tx_0)^{\oplus(n-m)\bv^0_k},
\end{multline*}
where $\tx_0\in\pi_0^{-1}(0)$. But decomposition~(\ref{seq0})
for the representation $\rho_0 \oplus \hrho(\tx_0)^{\oplus(n-m)\bv^0_k}$
gives
$$
\hrho(\tx_0)^{\oplus(n-m)\bv^0_k} \seq0 
\left(\bigoplus_{k=0}^{d-1}(\rho_k)\right)^{\oplus(n-m)\bv^0_k},
$$
hence $g_0(f_0(\rho))\seq0\rho$, the first equality 
of~(\ref{todo}) is proved.

Now, take $x_1,\dots,x_n\in X/\Gamma$. Then 
\begin{multline*}
f_0(g_0(x_1,\dots,x_n)) = 
f_0(\rho_0\oplus\hrho(x_1)\oplus\dots\oplus\hrho(x_n)) = \\
\supp(\Tot(\hrho(x_1))\oplus\dots\oplus\Tot(\hrho(x_n))) = 
\{x_1,\dots,x_n\}.
\end{multline*}
This completes the proof of the Proposition
\qed\medskip

\subsection*{$\Gamma$-equivariant Hilbert scheme}

The action of the group $\Gamma$ on $X$ induces an action 
of $\Gamma$ on the Hilbert scheme $X^{[nN]}$ of $nN$-tuples 
of points on $X$. Let $\left(X^{[nN]}\right)^\Gamma$ denote 
the $\Gamma$-invariant locus. Since the group $\Gamma$ acts 
on the tautological bundle $\CV$ on $X^{[nN]}$, the restriction
of $\CV$ to $\left(X^{[nN]}\right)^\Gamma$ decomposes into 
the direct sum of locally free sheaves indexed by the
set of irreducible representations of $\Gamma$:
\begin{equation}\label{udecomp}
\CV_{|\left(X^{[nN]}\right)^\Gamma} = 
\bigoplus_{k=0}^{d-1}\CV_k\otimes R_k.
\end{equation}
Let $X^{\Gamma[n]}$ denote the locus of $\left(X^{[nN]}\right)^\Gamma$,
where the rang of every $\CV_k$ equals~$n\bv^0_k$
(or, in the other words, where $\CV$ is a multiplicity of
the regular representation). It is clear that
$X^{\Gamma[n]}$ is a union of some of the connected components
of $\left(X^{[nN]}\right)^\Gamma$. We will reffer to $X^{\Gamma[n]}$
as the $\Gamma$-equivariant Hilbert scheme of $X$. 

Denote by $\cplus$ the positive octant of $\R^d$, that is
$$
\cplus = \{\zeta_\R\in\R^d\ |\ \zeta_\R^k>0\ \text{for all $0\le k\le d-1$}\}
$$
Note that any $\zeta$ with $\zeta_\R\in\cplus$ is generic (because 
its real component $\zeta_\R$ is stricktly positive) and any 
$\zeta_\R\in\cplus$ satisfies the restriction~(\ref{pos}).

The following Theorem is well known (see~\cite{W}, \cite{VV}).

\begin{theorem}\label{mplus}
For any integer $n\ge0$ and any $\zeta_\R\in\cplus$ the quiver 
variety $\FM_{(0,\zeta_R)}(n\cdot\bv^0,\bw^0)$ is isomorphic 
to the $\Gamma$-equivariant Hilbert scheme~$X^{\Gamma[n]}$.
\end{theorem}

We begin with the following evident Lemma.

\begin{lemma}\label{st}
Let $U$ be a $\Gamma$-equivariant finite dimensional $\C[X]$-module 
with a $\Gamma$-invariant vector $u\in U$ and let $\rho$ be 
the corresponding $(1,\bu)$-dimensional representation 
of the quiver $(Q(\bw^0),J(0))$. Then $\rho$ is 
$\chi_{\bu}(\zeta_\R)$-stable with $\zeta_\R\in\cplus$ 
iff $U$ has no proper $\Gamma$-equivariant $\C[X]$-submodule 
containing $u$.
\end{lemma}
{\sl Proof:}
Assume that $u\in U'\subset U$ is a $\Gamma$-equivariant proper
$\C[X]$-submodule. Let $\rho'$ be the corresponding 
$(1,\bu')$-dimensional representation of the quiver. 
Then $\rho'$ is a subrepresentation in $\rho$ and 
$$
\chi_\bu(\zeta_\R)(\rho') = 
(-\zeta_\R(\bu),\zeta_\R)(1,\bu') = 
\zeta_R(\bu') - \zeta_\R(\bu) = - \zeta_\R(\bu-\bu') < 0,
$$
since $\bu>\bu'$ and $\zeta_\R$ is positive. Hence $\rho$ is unstable. 

Similarly, assume that $\rho$ is unstable. Then it contains
a subrepresentation $\rho'$ such that
$$
\chi_\bu(\zeta_\R)(\dim\rho') =
(-\zeta_\R(\bu),\zeta_\R)(\bu'_*,\bu') = 
\zeta_R(\bu') - \bu'_*\zeta_\R(\bu) < 0,
$$
where $(\bu'_*,\bu') = \dim\rho'$.
Since $\zeta_\R$ is positive the case $\bu'_*=0$ is impossible.
Hence $\bu'_*=1$. This means that the subspace 
$\Tot(\rho')\subset\Tot(\rho) = U$  contains the vector $u$.
On the other hand, it is a proper $\Gamma$-equivariant
$\C[X]$-submodule in~$U$.
\qed\medskip

Now we can prove the Theorem.

{\sl Proof:}
Recall that $\CV$ carries the structure of a family of quotient algebras
of the algebra of functions $\C[X]$, and in particular, of 
$\Gamma$-equivariant $\C[X]$-module. In particular, the unit of
the algebra $\C[X]$ induces a $\Gamma$-equivariant morphism
$i:\CO\to\CV$. Applying relative analog of Corollary~\ref{vvk1} 
we see that for $\CV_*=\CO$ we obtain on $(\CV_*,\CV_k)$ the 
natural structure of a family of $(1,n\bv^0)$-dimensional representations 
of the quiver $(Q(\bw^0),J(0))$. 

Since $\CV$ is a family of quotient algebras of $\C[X]$
it follows that it is pointwise generated by the image of
$1\in\C[X]$ as a family of $\C[X]$-modules. In other words,
it has no proper $\Gamma$-invariant $\C[X]$-submodules
containig the image of~$1$. Hence Lemma~\ref{st} implies
$\chi_{n\bv^0}(\zeta_\R)$-stability of the family
$(\CV_*,\CV_k)$. Thus, we obtain a map
$$
f_+:X^{\Gamma[n]} \to 
\Mod_{Q(\bw^0),J(0)}((1,n\bv^0),\chi_{n\bv^0}(\zeta_\R)) = 
\FM_{(0,\zeta_\R)}(n\bv^0,\bw^0).
$$

Similarly, let $\urho = (\CV_*,\CV_k)$ be a universal family
on the moduli space 
$\Mod_{Q(\bw^0),J(0)}((1,n\bv^0),\chi_{n\bv^0}(\zeta_\R))$.
Twisting it by the line bundle $\CV^*_*$ we can 
assume that $\CV_*=\CO$. Then Corollary~\ref{vvk1} implies
that $\Tot(\urho)$ is a family of $\Gamma$-equivariant $\C[X]$-modules, 
endowed with a morphism $i:\CO\to\Tot(\urho)$. Furthermore, the 
$\chi_{n\bv^0}(\zeta_\R)$-stability of the family $\urho$ with 
Lemma~\ref{st} imply that $\Tot(\urho)$ is in fact a family of 
quotient algebras of $\C[X]$, when $1\in\C[X]$ goes to 
$i(1)\in\Tot(\urho)$. Since
$$
\rk(\Tot(\urho)) = \sum_{k=0}^{d-1}\rk(\CV_k)\dim R_k = 
\sum_{k=0}^{d-1}(n\bv^0_k)\bv^0_k = 
n\sum_{k=0}^{d-1}(\bv^0_k)^2 = nN
$$
we obtain a map
$$
g_+:\FM_{(0,\zeta_\R)}(n\bv^0,\bw^0) =
\Mod_{Q(\bw^0),J(0)}((1,n\bv^0),\chi_{n\bv^0}(\zeta_\R)) \to
X^{[nN]}.
$$
The $\Gamma$-equivariance of $\CV$ implies that $f_+$ 
goes in fact to the $\Gamma$-invariant part of $X^{[nN]}$
and moreover to $X^{\Gamma[n]}$.

It remains to notice that the maps $f_+$ and $g_+$
are evidently mutually inverse.
\qed\medskip

\subsection*{The Hilbert scheme of $X_\Gamma$}

Let $K_n,\BK_n\subset\Z^d$ denote the following finite sets 
of dimension vectors
$$
\begin{array}{c}
K_n = \{ \bv\in\Z_{\ge0}^d\ |\ \bv_0\bv^0 < \bv < n\bv^0 \}\smallskip\\
\BK_n = \{ \bv\in\Z_{\ge0}^d\ |\ \bv_0\bv^0 \le \bv \le n\bv^0 \}
\end{array}
$$
Then
$$
\cminus(n) = \left\{\zeta_\R\in\R^d\ \left|\ 
\begin{cases}
\zeta_\R(\bv^0)>0,\\
\zeta_\R(\bv) > n\zeta_\R(\bv^0),&\quad
\text{for all $\bv\in K_n$}
\end{cases}\right.\right\}
$$
is a convex polyhedral cone.

\begin{lemma}\label{cm}
We have
$$
\cminus(n) = \{\zeta_\R\in\R^d\ |\ 
\frac1n\min_{k=1}^{d-1}\zeta_\R^k-\sum_{k=1}^{d-1}\zeta_\R^k\bv^0_k > 
\zeta_\R^0 > - \sum_{k=1}^{d-1}\zeta_\R^k\bv^0_k \}
$$
In particular, the cone $\cminus(n)$ is nonempty and for any 
$\zeta_\R\in\cminus(n)$ we have $\zeta_\R^k>0$ for all $1\le k\le d-1$.
\end{lemma}
{\sl Proof:}
Easy.
\qed\medskip

\begin{theorem}\label{mminus}
For any $n\ge0$ and any $\zeta_\R\in\cminus(n)$ 
the quiver variety $\FM_{(0,\zeta_\R)}(n\bv^0,\bw^0)$ is isomorphic
to the Hilbert scheme $X_\Gamma^{[n]}$.
\end{theorem}
{\sl Proof:}
Let $\urho=\urho^{\zeta_\R}=(\CV_*,\CV_k,\uB_h,\ui_0,\uj_0)$ 
be the universal family on $X_\Gamma$ normalized as in Lemma~\ref{norm}. 
Let $Z\subset X_\Gamma$ be a length $n$ subscheme in $X_\Gamma$ and 
let $\CO_Z$ be its structure sheaf. We associate to $Z$ the following 
representation $\rho(Z)=(V_*(Z),V_k(Z),B(Z),i_0(Z),j_0(Z))$ 
of the quiver $Q(\bw^0)$. We put
$$
\begin{array}{l}
V_k(Z) = \Gamma(X_\Gamma,\CV_k\otimes\CO_Z)\qquad k=0,\dots,d-1,\\
V_*(Z)= \Gamma(X_\Gamma,\CO_{X_\Gamma}) \cong \C
\end{array}
$$
with morphisms $B_h(Z):V_{out(h)}(Z)\to V_{in(h)}(Z)$ induced by the maps
$\uB_h\otimes\id_{\CO_Z}:\CV_{out(h)}\otimes\CO_Z\to\CV_{in(h)}\otimes\CO_Z$, 
while $j_0(Z)=0$ and $i_0(Z)$ is induced by the composition
$$
\CO_{X_\Gamma} \cong \CV_0 \to \CV_0\otimes\CO_Z,
$$
where the first is the isomorphism of Lemma~\ref{norm}, 
while the second morphism is obtained from the canonical projection 
$\CO_{X_\Gamma}\to\CO_Z$ by tensoring with~$\CV_0$.

Since $Z$ is a length $n$ subscheme of $X_\Gamma$ and $\rk\CV_k=\bv^0_k$ 
it follows that $V_k(Z)$ is an $n\bv^0_k$-dimensional vector space
for any $0\le k\le d-1$. On the other hand, $V_*(Z)$ is 
$1$-dimensional by definition. Hence $\rho(Z)$ is 
$(1,n\bv^0)$-dimensional reresentation of the quiver $Q(\bw^0)$.
Further, it follows from~(\ref{mmc}) and Lemma~\ref{li0j0} that
$\rho(Z)$ satisfies the relations $J(0)$.
Now we want to check that $\rho(Z)$ is
$\chi_{n\bv^0}(\zeta_\R)$-stable for any $\zeta_\R\in\cminus(n)$.
We begin with the following Proposition.

\begin{proposition}\label{subrep}
If $\rho'\subset\rho(Z)$ is a subrepresentation and
$\rho'\ne\rho(Z)$ then we have
$$
\dim\rho' = (0,\bv)\text{ with $\bv\in\BK_n$, or }
\dim\rho' = (1,\bv)\text{ with $\bv\in  K_n$.}
$$
\end{proposition}
{\sl Proof:}
Choose a chain of subschemes 
$$
Z = Z_n\supset Z_{n-1}\supset \dots \supset Z_1\supset Z_0=\emptyset,
$$
such that $Z_i$ is a length $i$ subscheme. This chain of subschemes
induces a chain of surjections of structure sheaves 
$$
\CO_Z = \CO_{Z_n}\to \CO_{Z_{n-1}}\to \dots 
\to \CO_{Z_1}\to \CO_{Z_0} = 0,
$$
and of representations of quivers
\begin{equation}\label{rhosur}
\rho(Z) = \rho(Z_n) \to \rho(Z_{n-1}) \to \dots 
\to \rho(Z_1) \to \rho(Z_0) = \rho_0.
\end{equation}
Note that for any $i=1,\dots,n$ we have exact sequence of sheaves
$$
0 \to \CO_{x_i} \to \CO_{Z_i} \to \CO_{Z_{i-1}} \to 0,
$$
for some $x_i\in X_\Gamma$, which induces exact sequence of representations
$$
0 \to \hrho(x_i) \to \rho(Z_i) \to \rho(Z_{i-1}) \to 0.
$$

Now assume that $\rho'$ is a subrepresentation in $\rho(Z)$.
Let $\rho'_i$ denote the image of $\rho'$ in $\rho(Z_i)$ with
respect to the surjection~(\ref{rhosur}). Then~(\ref{rhosur})
induces a chain of surjections
$$
\rho' = \rho'_n \to \rho'_{n-1} \to \dots \to \rho'_1 \to \rho'_0 = \rho_0.
$$
Let $\rho''_i$ denote the kernel of the map $\rho'_i\to\rho'_{i-1}$.
Then $\rho''_i$ is a subrepresentation in $\hrho(x_i)$, hence 
$\dim\rho''_i=(0,\bu^i)$ for some dimension vector $\bu^i\le\bv^0$.
Since $\bv^0_0=1$, it follows that either $\bu^i_0=0$ or $\bu^i_0=1$.

Assume that $\bu^i_0=1$. Then it is easy to see that we can extend 
the subrepresentation $\rho''_i\subset\hrho(x_i)$
to $(1,\bu^i)$-dimensional subrepresentation 
$\widetilde{\rho''_i}\subset\rho(x_i)$. Now, note that
if $\bu^i<\bv^0$ then since $\bu^i_0=\bv^0_0$ and 
since $\zeta_\R^k>0$ for all $1\le k\le d-1$ by Lemma~\ref{cm},
it follows that
$$
\chi_{\bv^0}(\zeta_\R)(1,\bu^i) = 
(-\zeta_\R(\bv^0),\zeta_\R)(1,\bu^i) <
(-\zeta_\R(\bv^0),\zeta_\R)(1,\bv^0) = 0,
$$
which contradicts the $\chi_{\bv^0}(\zeta_\R)$-stability of $\rho(x_i)$.
Thus we have proved that if $\bu^i_0=1$ then $\bu^i=\bv^0$.

On the other hand, if $\bu^i_0=0$ then certainly $\bu^i\ge 0$. 
Thus in both cases we have
$$
\bu^i \ge \bu^i_0\bv^0.
$$
Summing up these inequalities from $i=1$ to $n$ we see that
$$
\bv = \sum_{i=1}^n\bu^i \ge \sum_{i=1}^n\bu^i_0\bv^0 = \bv_0\bv^0,
$$
hence $\bv\in\BK_n$.

It remains to check that if $\dim\rho'=(1,\bv)$ then $\bv\in K_n$. 
So, assume that it is not true. Then $\bv=m\bv^0$ for some $0\le m\le n$. 
Let us check that the case $m<n$ is impossible. 

For this we will use induction in $n$. The base of induction, 
$n=0$, is trivial. So assume that $n>0$, $Z\subset X_\Gamma$
is a length $n$ subscheme, and that $\rho'\subset\rho(Z)$ is 
a $(1,m\bv^0)$-dimensional subrepresentation with $m<n$. 
Consider the following commutative diagram with exact rows
$$
\begin{CD}
0 @>>>  \rho''_n   @>>> \rho'   @>>>    \rho'_{n-1}   @>>>      0       \\
@.      @VVV            @VVV            @VVV                            \\
0 @>>>  \hrho(x_n) @>>> \rho(Z) @>>>    \rho(Z_{n-1}) @>>>      0       
\end{CD}
$$
where $\rho''_n$ and $\rho'_{n-1}$ was defined above. Denote
$$
\dim\rho''_n = (0,\bu),\qquad
\dim\rho'_{n-1} = (1,\bv').
$$
Then the above arguments show that
$$
\bu  \ge \bu_0\bv^0,\qquad
\bv' \ge \bv'_0\bv^0.
$$
But since $\bu + \bv' = m\bv^0 = (\bu_0+\bv'_0)\bv^0$ it follows that
$$
\bu = m''\bv^0,\qquad
\bv'= m' \bv^0,\qquad m'+m''=m.
$$
In particular, $\rho(Z_{n-1})$ contains a $(1,m'\bv^0)$-dimensional
subrepresentation. The induction hypothesis for $Z_{n-1}$ then
implies that $m'=n-1$. On the other hand, we have $m''=1$ or $m''=0$.
In the first case we have $m=n$, a contradiction. In the second
case, it follows that $\rho'=\rho(Z_{n-1})$, hence the exact sequence
$$
0 \to \hrho(x_n) \to \rho(Z) \to \rho(Z_{n-1}) \to 0
$$
splits. Comparing the definition of $\rho(Z)$ with the definition 
of the equivalence $\Psi_{\zeta_\R}$, we see that the splitting
of the above sequence implies that the sequence
$$
0 \to \CO_{x_n} \to \CO_{Z} \to \CO_{Z_{n-1}} \to 0,
$$
also splits, and that this splitting is compatible with the projections
$\CO_{X_\Gamma}\to\CO_Z$ and $\CO_{X_\Gamma}\to\CO_{Z_{n-1}}$.
But this means that
$$
J_{Z_{n-1}} \cong J_Z \oplus \CO_{x_n},
$$
which is false. Thus we again come to a contradiction,
and the Lemma is proved.
\qed\medskip

Now we can finish the proof of the Theorem. 
The $\chi_{n\bv^0}(\zeta_\R)$-stability of 
the representation $\rho(Z)$ follows immediately 
form the definition of the cone $\cminus(n)$ and
from Proposition~\ref{subrep}. Thus $\rho(Z)$ form
a family of $(1,n\bv^0)$-dimensional $\chi_{n\bv^0}(\zeta_\R)$-stable
representations of the quiver $(Q(\bw^0),J(0))$ 
over the Hilbert scheme $X_\Gamma^{[n]}$. This family induces a map
$$
g_-:X_\Gamma^{[n]} \to 
\Mod_{Q(\bw^0),J(0)}((1,n\bv^0),\chi_{n\bv^0}(\zeta_\R)) =
\FM_{(0,\zeta_\R)}(n\bv^0,\bw^0).
$$
Moreover, it is easy to see that the map fits into the following
commutative diagram:
$$
\begin{CD}
X_\Gamma^{[n]} @>{g_-}>>        \FM_{(0,\zeta_\R)}(n\bv^0,\bw^0) \\
@VVV                            @VV{\pi_0}V                      \\
S^n(X/\Gamma)  @>{g_0}>>        \FM_{(0,0)}(n\bv^0,\bw^0)
\end{CD}
$$
where the left vertical arrow is the composition of the 
Hilbert-Chow morphism $X_\Gamma^{[n]}\to S^n(X_\Gamma)$ with 
the map induced by the projection $\pi_0:X_\Gamma\to X/\Gamma$.
Since both vertical maps are birational, and $g_0$ is
an isomorphism by Proposition~\ref{psi0} it follows that
$g_-$ is a regular birational map. On the other hand,
both $X_\Gamma^{[n]}$ and $\FM_{(0,\zeta_\R)}(n\bv^0,\bw^0)$ are 
holomorphically symplectic varieties, hence any regular birational 
map between them is an isomorphism.
\qed\medskip

\begin{remark}\label{pn}
H.~Nakaijma indicated to me the direction of his arguments.
He constructs a map $\FM_{(0,\zeta_\R)}(n\bv^0,\bw^0)\to X_\Gamma^{[n]}$
(note that the direction is opposite to that of the map $g_-$)
using certain complex and utilizing the stability condition 
to ensure that certain cohomology groups vanish.
\end{remark}

Combining Theorems~\ref{mplus} and \ref{mminus} with Theorem~\ref{na}
we get the following Corollary.

\begin{corollary}\label{cor}
The $\Gamma$-equivariant Hilbert scheme $X^{\Gamma[n]}$ is 
diffeomorphic to the Hilbert scheme $X_\Gamma^{[n]}$. 
In particular, we have an isomorphism of cohomology groups
$$
H^\bullet(X^{\Gamma[n]},\Z) \cong H^\bullet(X_\Gamma^{[n]},\Z).
$$
\end{corollary}

\begin{remark}
It is also easy to deduce that for any 
$$
\zeta_\R \in \{\zeta_\R\ |\ \zeta_\R(\bv^0) = 0,\quad\text{and}\quad
\zeta_\R^k > 0\quad\text{for all $1\le k\le d-1$}\}
$$
we have
$$
\FM_{(0,\zeta_\R)}(n\bv^0,\bw^0) \cong S^n(X_\Gamma).
$$
\end{remark}

\subsection*{Generalizations: the Calogero-Moser space}

Let $\Gamma$ be a finite subgroup in $SL(2,\C)$ and $\Gamma'\subset\Gamma$
its central subgroup. Let $\FM^\Gamma$ and $\FM^{\Gamma'}$ denote 
the quiver varieties corresponding to the affine Dynkin graphs
of $\Gamma$ and $\Gamma'$ respectively. 

Since $\Gamma'$ is central in $\Gamma$ the embedding $\Gamma'\to\Gamma$
induces embeddings 
$$
\sigma_\C^*:Z(\C[\Gamma']) \to Z(\C[\Gamma]),\qquad
\sigma_\R^*:Z(\R[\Gamma']) \to Z(\R[\Gamma]).
$$
\begin{remark}
In terms of the Dynkin graphs the embedding $\sigma^*$ can be described
as follows. Let $I$ and $I'$ be the sets (of isomorphism classes) of 
irreducible representations of $\Gamma$ and $\Gamma'$ respectively.
In other words, $I$ and $I'$ are the sets of vertices of the corresponding
affine Dynkin graphs. Consider an irreducible $\Gamma$-module $R_i^\Gamma$,
$i\in I$. Since $\Gamma'$ is central it follows from the Schur Lemma that 
the restriction $(R_i^\Gamma)_{|\Gamma'}$ is a multiplicity of an irreducible 
representation of the group $\Gamma'$, say~$R^{\Gamma'}_{i'}$, $i'\in I'$. 
So, associating this way to arbitrary vertex $i\in I$ the vertex $i'\in I'$
we obtain a map $\sigma:I\to I'$. The map $\sigma^*$ is induced by~$\sigma$.
\end{remark}

Let $N=|\Gamma/\Gamma'|$ denote the index of $\Gamma'$ in $\Gamma$. Choose 
arbitrary generic $\zeta\in Z(\C[\Gamma'])\oplus Z(\R[\Gamma'])$.
Let also $W$ be an arbitrary representation of $\Gamma$, let 
$W'$ be the restriction of $W$ to $\Gamma'$ and let $V'$ be
an arbitrary representation of $\Gamma'$. Let $\bv'$, $\bw'$
and $\bw$ denote the classes of $V'$, $W'$ and $W$ in the Grothendieck
ring of $\Gamma'$ and $\Gamma$ respectively (i.e.~their dimension vectors).

Recall that by Lemma~\ref{gqv} the quiver variety 
$\FM=\FM^{\Gamma'}_\zeta(\bv',\bw')$ coincides with the set of all triples
$$
(B,i,j)\in\Hom_{\Gamma'}(V'\otimes L,V')\oplus
\Hom_{\Gamma'}(W',V')\oplus\Hom_{\Gamma'}(V',W')
$$
such that 
$$
[B,B]+ij=-\zeta_\C, \qquad [B,B^\dagger]+ii^\dagger-j^\dagger j=-\zeta_\R
$$
and modulo action of $U_{\Gamma'}(V')$.
We define for every $\gamma\in\Gamma$ another triple 
$$
B^\gamma = B \cdot (1\otimes\gamma) : 
V\otimes L @>{1\otimes\gamma}>> V\otimes L @>{B}>> V,\quad
i^\gamma = i\cdot\gamma,\quad j^\gamma = \gamma^{-1}\cdot j.
$$
Since $\Gamma$ commutes with $\Gamma'$, it follows 
that $(B^\gamma,i^\gamma,j^\gamma)$ gives another point of $\FM$,
hence the correspondence $B\mapsto B^\gamma$ defines an action 
of the group $\Gamma$ on the quiver variety~$\FM$. Let $\FM^\Gamma$
denote the fixed points set of $\Gamma$ on $\FM$.

Take arbitrary $(B,i,j)\in\FM^\Gamma$. Then $(B,i,j)$ and
$(B^\gamma,i^\gamma,j^\gamma)$ should be conjugated under 
the action of $U_{\Gamma'}(V')$. Hence there exists 
$g_\gamma\in U_{\Gamma'}(V')$ such that 
\begin{equation}\label{gg}
g_\gamma B g_\gamma^{-1} = B^\gamma,\quad
g_\gamma i = i^\gamma,\quad
j g_\gamma^{-1} = j^\gamma.
\end{equation}
Moreover, when $\zeta$ is generic such $g_\gamma$ is unique
(because the action of $U_{\Gamma'}(V')$ is free in this case).
It follows that $\gamma\mapsto g_\gamma$ defines an action 
of $\Gamma$ on $V'$ extending the action of $\Gamma'$.

Now let $V$ be an arbitrary representation of $\Gamma$,
such that its restriction to $\Gamma'$ is isomorphic to $V'$.
Let $\bv$ be its class in the Grothendieck ring of $\Gamma$.
Let $\FM^\Gamma_\bv$ denote the locus of the set $\FM^\Gamma$,
where the defined above structure of a representation of $\Gamma$
on $V'$ is isomorphic to $V$. Then we have.

\begin{theorem}\label{gen}
For any generic $\zeta\in Z(\C[\Gamma'])\oplus Z(\R[\Gamma'])$, 
$\bv,\bw\in K_0(\Gamma)$ let 
$\bv'=\bv_{|\Gamma'}$, $\bw'=\bw_{|\Gamma'}$. Then we have
$$
\left(\FM^{\Gamma'}_\zeta(\bv',\bw')\right)^\Gamma_\bv =
\FM^\Gamma_{\sigma^*\zeta}(\bv,\bw).
$$
\end{theorem}
{\sl Proof:}
Note that the equations~(\ref{gg}) mean that the triple $(B,i,j)$
is $\Gamma$-equivariant with respect to the action $\gamma\to g_\gamma$
of the group $\Gamma$ on $V$ and its canonical actions on $W$ and $L$.
Choose an arbitrary $\Gamma$-equivariant isomorphism $V'\to V$.
Then 
$$
(B,i,j)\in\Hom_\Gamma(V\otimes L,V)\oplus
\Hom_\Gamma(W,V)\oplus\Hom_\Gamma(V,W).
$$
Furthermore, it follows from the definition of $\sigma^*$ that
$$
[B,B]+ij=-\sigma_\C^*\zeta_\C, \qquad 
[B,B^\dagger]+ii^\dagger-j^\dagger j=-\sigma_\R^*\zeta_\R.
$$
Thus we obtain a map 
$$
\left(\FM^{\Gamma'}_\zeta(\bv',\bw')\right)^\Gamma_\bv \to
\FM^\Gamma_{\sigma^*\zeta}(\bv,\bw).
$$
Similarly, for any $(B,i,j)\in\FM^\Gamma_{\sigma^*\zeta}(\bv,\bw)$
forgetting the $\Gamma$-structure on $V$ and $W$ we can consider it
as a point of $\FM^{\Gamma'}_\zeta(\bv',\bw')$. Thus we obtain
the inverse map.
\qed\medskip

Consider the case $\Gamma'=\{1\}$, $W=R_0$, $W'=\C$, $V'=\C^{nN}$. 
Then the quiver variety $\FM^{\Gamma'}_\zeta(\bv',\bw')$ 
coincides with the Hilbert scheme $X^{[nN]}$ when $\zeta_\C=0$ 
and with the so-called Calogero-Moser space $\CM_{nN}$ when 
$\zeta_\C\ne0$. Thus in the case $\zeta_\C=0$, Theorem~\ref{gen} 
specializes to Theorem~\ref{mplus}, and in the case $\zeta_\C=\tau\ne0$ 
we obtain the following.

\begin{corollary}
For any $\tau\ne0$ and $n\ge0$ we have
$$
(\CM_{nN})^\Gamma_{n\bv^0} = \FM^\Gamma_\zeta(n\bv^0,\bw^0),
$$
where $\zeta_\C=(\tau,\dots,\tau)$ and $\zeta_\R$ is arbitrary.
\end{corollary}

\section{Combinatorial applications}

From now on assume that $\Gamma\cong\Z/d\Z$ 
(the $\widetilde{A}_{d-1}$-case). Then both $X$ and $X_\Gamma$
admit a $\Gamma$-equivariant action of the torus $\C^*\times\C^*$.
The first action is a coordinatewise dilation, and the second
is induced by the first one. These actions induce actions on
the Hilbert schemes $X^{\Gamma[n]}$ and $X_\Gamma^{[n]}$.
In both cases, there are only finite number of fixed points.
Now we will give their combinatorial description.

Let us consider a Young diagram as a domain in the top right-hand 
octant of the coordinate plane. We associate to a box of a diagram 
the coordinates of its bottom left-hand corner. Thus the coordinates 
of any box are nonnegative integers. We denote by $(p,q)$ the box with 
the coordinates $(p,q)$. A Young diagram $\Delta$ is called {\em uniformly
coloured in~$d$ colours} if the integer
$$
n_i(\Delta) = \#\{(p,q)\in\Delta\ |\ p-q\equiv i\pmod{d}\},
\qquad0\le i\le d-1
$$
(the number of boxes of the colour $i$) doesn't depend on $i$. One can 
say that we colour the diagonals of the digram $\Delta$ periodically 
in $d$ colours and call the diagram uniformly coloured if it has equal 
number of boxes of each colour. Let $UCY(n,d)$ denote the number of 
uniformly coloured in $d$ colours Young diagrams with $n$ boxes of 
each colour.

\begin{lemma}
The $\C^*\times\C^*$-fixed points on the $\Gamma$-equivariant
Hilbert scheme $X^{\Gamma[n]}$ are in a $1-1$ correspondence 
with uniformely coloured in $d$ colours Young diagrams 
with $n$ boxes of each colour.
\end{lemma}
{\sl Proof:}
Note that by definition $X^{\Gamma[n]}$ is a connected component 
of $(X^{[dn]})^\Gamma$. The fixed points on $X^{[dn]}$ are numbered
by Young diagrams with $dn$ boxes. All these points are $\Gamma$-invariant
(because the action of $\Gamma$ factors through the torus action),
so it remains to understand, which of these points lie in the
component. 

Note that the fiber of the tautological bundle over $X^{[dn]}$ 
at the $\C^*\times\C^*$-fixed point corresponding to a Young 
diagram~$\Delta$ is isomorphic as a $\Gamma$-module to the
representation $\oplus_{(p,q)\in\Delta} R_{p-q}$. In particular,
it is a multiplicity of the regular representation iff $\Delta$
is uniformly coloured in $d$ colours.
\qed\medskip

Let $CY(n,d)$ be the number of all ordered
collections $(\Delta_1,\dots,\Delta_d)$ of Young diagrams such that
$$
\sum_{k=1}^d |\Delta_k|=n,
$$
where $|\Delta|$ is the number of boxes in $\Delta$. 

\begin{lemma}
The $\C^*\times\C^*$-fixed points on the Hilbert scheme 
$X_\Gamma^{[n]}$ are in a $1-1$ correspondence 
with ordered collections $(\Delta_1,\dots,\Delta_d)$ of
Young diagrams, such that $\sum_{k=1}^d|\Delta_k|=n$.
\end{lemma}
{\sl Proof:}
First note that the number of $\C^*\times\C^*$ fixed
points on~$X_\Gamma$ equals~$d$. It follows for example, 
from the evident fact, that $X_\Gamma=X^{\Gamma,1}$, and 
$UCY(1,d)=d$. Further, denoting these points by $x_1,\dots, x_d$ 
it is easy to see that $Z$ is a fixed point on $X_\Gamma^{[n]}$ 
iff it splits as a union $Z=Z_1\cup\dots\cup Z_d$, where $Z_i$ is 
a $\C^*\times\C^*$-invariant length $n_i$-subscheme in $X_\Gamma$
with support at $x_i$. Linearizing the action of the torus in
a neighbourhood of $x_i\in X_\Gamma$ we see that such $Z_i$ are
numbered by Young diagrams $\Delta_i$ with $|\Delta_i|=n_i$.
Finally, the condition $Z\in X_\Gamma^{[n]}$ is just $\sum_{i=1}^dn_i=n$,
hence the Lemma.
\qed\medskip

\begin{theorem}
For any $n,d > 0$ we have $UCY(n,d) = CY(n,d)$.
\end{theorem}
\begin{remark}
See \cite{JK} for the combinatorial proof of this identity.
\end{remark}
{\sl Proof:}
Using the arguments of Nakajima (see \cite{Na2}, Chapter 5)
one can check that the dimension of the cohomolgy groups
of the Hilbert schemes $X^{\Gamma[n]}$ and $X_\Gamma^{[n]}$
are equal to the number of the fixed points with respect to
the torus action. Thus
$$
H^\bullet(X^{\Gamma[n]},\Z) \cong \Z^{UCY(n,d)},\qquad
H^\bullet(X_\Gamma^{[n]},\Z) \cong \Z^{CY(n,d)}
$$
and the Theorem follows from Corollary~\ref{cor}.
\qed

\end{document}